\documentclass[12pt]{article}      
\usepackage{amssymb}
\usepackage{eucal}
\usepackage{amsmath}
\usepackage{amscd}

\newtheorem{theorem}{Theorem}[section]
\newtheorem{prop}[theorem]{Proposition}
\newtheorem{defn}[theorem]{Definition}
\newtheorem{lemma}[theorem]{Lemma}
\newtheorem{coro}[theorem]{Corollary}
\newtheorem{prop-def}{Proposition-Definition}[section]

\newtheorem{remark}{Remark}[section]

\newcommand{\nc}{\newcommand}

\nc{\bin}[2]{ (_{\stackrel{\scs{#1}}{\scs{#2}}})}  
\nc{\binc}[2]{ \left (\!\! \begin{array}{c} \scs{#1}\\
    \scs{#2} \end{array}\!\! \right )}  
\nc{\bincc}[2]{  \left ( {\scs{#1} \atop
    \vspace{-1cm}\scs{#2}} \right )}  
\nc{\bs}{\bar{S}}
\nc{\la}{\longrightarrow}
\nc{\rar}{\rightarrow}
\nc{\dar}{\downarrow}
\nc{\dap}[1]{\downarrow \rlap{$\scriptstyle{#1}$}}
\nc{\uap}[1]{\uparrow \rlap{$\scriptstyle{#1}$}}
\nc{\defeq}{\stackrel{\rm def}{=}}
\nc{\dis}[1]{\displaystyle{#1}}
\nc{\dotcup}{\ \displaystyle{\bigcup^\bullet}\ }
\nc{\hcm}{\ \hat{,}\ }
\nc{\hts}{\hat{\otimes}}
\nc{\hcirc}{\hat{\circ}}
\nc{\lleft}{\left [ \begin{array}{c} {} \\ {} \end{array}
    \right .  \!\!\!\!\!\!\!} 
\nc{\lright}{ \!\!\!\!\!\!\!
    \left . \begin{array}{c} {} \\ {} \end{array}
    \right ] }
\nc{\curlyl}{\left \{ \begin{array}{c} {} \\ {} \end{array}
    \right .  \!\!\!\!\!\!\!} 
\nc{\curlyr}{ \!\!\!\!\!\!\!
    \left . \begin{array}{c} {} \\ {} \end{array}
    \right \} }
\nc{\longmid}{\left | \begin{array}{c} {} \\ {} \end{array}
    \right . \!\!\!\!\!\!\!}
\nc{\ora}[1]{\stackrel{#1}{\rar}}
\nc{\ola}[1]{\stackrel{#1}{\la}}
\nc{\scs}[1]{\scriptstyle{#1}}
\nc{\sss}{\subsubsection}
\nc{\mrm}[1]{{\rm #1}}
\nc{\margin}[1]{\marginpar{\rm #1}}   
\nc{\dirlim}{\displaystyle{\lim_{\longrightarrow}}\,}
\nc{\invlim}{\displaystyle{\lim_{\longleftarrow}}\,}
\nc{\mvp}{\vspace{0.3cm}}
\nc{\tk}{^{(k)}}
\nc{\tp}{^\prime}
\nc{\ttp}{^{\prime\prime}}
\nc{\svp}{\vspace{2cm}}
\nc{\vp}{\vspace{8cm}}
\nc{\proofbegin}{\noindent{\bf Proof: }}
\nc{\proofend}{$\blacksquare$ \vspace{0.3cm}}
\nc{\modg}[1]{\!<\!\!{#1}\!\!>}
\nc{\intg}[1]{F_C(#1)}
\nc{\lmodg}{\!<\!\!}
\nc{\rmodg}{\!\!>\!}
\nc{\cpi}{\widehat{\Pi}}
\nc{\sha}{{\mbox{\cyr X}}}  
\nc{\shpr}{\diamond}    
\nc{\vep}{\varepsilon}
\nc{\labs}{\mid\!}
\nc{\rabs}{\!\mid}
\nc{\lpair}[2]{{<{#1} \mid {#2}>_\lambda\ }}
\nc{\Cq}{{C\!\!<\!\! \frakq\!\!>}}

\nc{\ann}{\mrm{ann}}
\nc{\Aut}{\mrm{Aut}}
\nc{\can}{\mrm{can}}
\nc{\colim}{\mrm{colim}}
\nc{\Cont}{\mrm{Cont}}
\nc{\rchar}{\mrm{char}}
\nc{\cok}{\mrm{coker}}
\nc{\dtf}{{R-{\rm tf}}}
\nc{\dtor}{{R-{\rm tor}}}

\nc{\Div}{{\mrm Div}}
\nc{\End}{\mrm{End}}
\nc{\Ext}{\mrm{Ext}}
\nc{\Fil}{\mrm{Fil}}
\nc{\Frob}{\mrm{Frob}}
\nc{\Gal}{\mrm{Gal}}
\nc{\GL}{\mrm{GL}}
\nc{\Hom}{\mrm{Hom}}
\nc{\hsr}{\mrm{H}}
\nc{\hpol}{\mrm{HP}}
\nc{\id}{\mrm{id}}
\nc{\im}{\mrm{im}}
\nc{\incl}{\mrm{incl}}
\nc{\length}{\mrm{length}}
\nc{\mchar}{\rm char\ }
\nc{\mpart}{\mrm{part}}
\nc{\ql}{{\QQ_\ell}}
\nc{\qp}{{\QQ_p}}
\nc{\rank}{\mrm{rank}}
\nc{\rcot}{\mrm{cot}}
\nc{\rdef}{\mrm{def}}
\nc{\rdiv}{{\rm div}}
\nc{\rtf}{{\rm tf}}
\nc{\rtor}{{\rm tor}}
\nc{\res}{\mrm{res}}
\nc{\SL}{\mrm{SL}}
\nc{\Spec}{\mrm{Spec}}
\nc{\tor}{\mrm{tor}}
\nc{\Tr}{\mrm{Tr}}
\nc{\tr}{\mrm{tr}}

\nc{\ab}{\mathbf{Ab}}
\nc{\Alg}{\mathbf{Alg}}
\nc{\Bax}{\mathbf{Bax}}
\nc{\bfk}{{\bf k}}
\nc{\bfone}{{\bf 1}}
\nc{\base}[1]{{a_{#1}}}
\nc{\detail}{\marginpar{\bf More detail}
    \noindent{\bf Need more detail!}
    \svp}
\nc{\Diff}{\mathbf{Diff}}   
\nc{\gap}{\marginpar{\bf Incomplete}\noindent{\bf Incomplete!!}
    \svp}
\nc{\FMod}{\mathbf{FMod}}
\nc{\Int}{\mathbf{Int}}
\nc{\Mon}{\mathbf{Mon}}
\nc{\remarks}{\noindent{\bf Remarks: }}
\nc{\Rep}{\mathbf{Rep}}
\nc{\Rings}{\mathbf{Rings}}
\nc{\Sets}{\mathbf{Sets}}
\nc{\bill}[1]{\marginpar{\bf To Bill}\noindent{\bf To Bill:}
    {\tt #1}\\ }
\nc{\li}[1]{\marginpar{\bf To Li}\noindent{\bf To Li:}
    {\tt #1}\\ }
\nc{\mlabel}[1]{\label{#1}}  
\nc{\dfootnote}[1]{{}}          

\nc{\BA}{{\mathbb A}}
\nc{\CC}{{\mathbb C}}
\nc{\DD}{{\mathbb D}}
\nc{\EE}{{\mathbb E}}
\nc{\FF}{{\mathbb F}}
\nc{\GG}{{\mathbb G}}
\nc{\HH}{{\mathbb H}}
\nc{\LL}{{\mathbb L}}
\nc{\NN}{{\mathbb N}}
\nc{\QQ}{{\mathbb Q}}
\nc{\RR}{{\mathbb R}}
\nc{\TT}{{\mathbb T}}
\nc{\VV}{{\mathbb V}}
\nc{\ZZ}{{\mathbb Z}}

\nc{\cala}{{\cal A}}
\nc{\calc}{{\cal C}}
\nc{\cald}{\mathcal{D}}
\nc{\cale}{{\cal E}}
\nc{\calf}{{\cal F}}
\nc{\calg}{{\cal G}}
\nc{\calh}{{\cal H}}
\nc{\cali}{{\cal I}}
\nc{\call}{{\cal L}}
\nc{\calm}{{\cal M}}
\nc{\caln}{{\cal N}}
\nc{\calo}{{\cal O}}
\nc{\calp}{{\cal P}}
\nc{\calr}{{\cal R}}
\nc{\calt}{{\cal T}}
\nc{\calw}{{\cal W}}
\nc{\calx}{{\cal X}}
\nc{\CA}{\mathcal{A}}

\nc{\fraka}{{\mathfrak a}}
\nc{\frakB}{{\mathfrak B}}
\nc{\frakm}{{\mathfrak m}}
\nc{\frakp}{{\mathfrak p}}
\nc{\frakq}{{\mathfrak q}}
\nc{\frakQ}{{\mathfrak Q}}

\font\cyr=wncyr10

\title{ 
Baxter Algebras and the Umbral Calculus
\thanks{MSC: 16W99, 05A50}}

\author{Li Guo\\
Department of Mathematics and Computer Science\\
Rutgers University at Newark}
\date{}
\begin{document}

\maketitle

\begin{abstract}
We apply Baxter algebras to the study of the umbral calculus. 
We give a characterization of the umbral calculus in terms 
of Baxter algebra. 
This characterization leads to a natural generalization 
of the umbral calculus that include the classical umbral calculus 
in a family of $\lambda$-umbral calculi parameterized 
by $\lambda$ in the base ring. 
\end{abstract}

\setcounter{section}{0}

\section{Introduction}
\mlabel{intro}

Baxter algebra and the umbral calculus are two areas 
in mathematics 
that have interested Rota throughout his life time 
and in which he has made 
prominent contributions. We will show in this paper 
that these two areas are intimately related.

The umbral calculus is the study and application of 
sequences of polynomials of binomial type and other 
related sequences. 
More precisely, 
a sequence $\{p_n (x)\mid n\in\NN \}$ of polynomials in 
$C[x]$ is called { a sequence of binomial type} if 
\[ p_n(x+y)= \sum_{k=0}^n p_k(x) p_{n-k}(y), 
\forall y\in C, n\in \NN.\]
Such sequences have fascinated mathematicians since 
the 19th century and include many of the most well-known 
sequences, such as the ones named after  
Abel, Bernoulli, Euler, Hermite and Mittag-Leffler.  
Even though polynomials of binomial type proved to be 
useful in several areas of mathematics, the foundation 
of umbral calculus was not firmly established for over 
a century since its first introduction. 
This situation changed completely in 1964 when 
G.-C. Rota\cite{Ro1} 
indicated that the theory can be rigorously formulated 
in terms of the algebra 
of functionals defined on the polynomials, later known as 
the umbral algebra. Rota's pioneer work was completed 
in the next decade by Rota and his collaborators 
\cite{RKO,RR,Rom}. 
Since then, there have been a number of generalizations 
of the umbral calculus~\cite{Rom,Ch,Lo,Me,Ve}. 

During the same period of time in which 
Rota embarked on laying down
the foundation of the umbral calculus, he also started the 
algebraic study of Baxter algebra which was first introduced 
by Baxter in connection with fluctuation theory in 
probability~\cite{Ba}. 
Fundamental to the study of Baxter algebra 
are the important works of Rota \cite{Ro2} and 
Cartier~\cite{Ca} that gave constructions of free Baxter 
algebras. By using a generalization of shuffle product 
in topology and geometry, the present author and W. Keigher 
gave another construction of free Baxter algebras
\cite{G-K1,G-K2}. This construction is applied to the 
further study of free Baxter algebras \cite{Gu1,Gu2,AGKO}. 

Our first purpose in this paper is to give a characterization 
of umbral calculus in terms of free Baxter algebras. 
We show that the umbral algebra is the 
free Baxter algebra of weight zero on the empty set.  
We also characterize the polynomial sequences studied in 
the umbral calculus in terms of operations in free Baxter 
algebras.

The second purpose of this paper is to use the free Baxter 
algebra formulation of the umbral calculus we have obtained 
to give a 
generalization of the umbral calculus, the 
$\lambda$-umbral 
calculus, for each constant $\lambda$ in the base ring $C$. 
The umbral calculus 
of Rota is the special case when $\lambda=0$. 

For simplicity, we only consider sequences of binomial 
type in this paper. The study of other sequences in 
the umbral calculus, 
such as Sheffer sequences and cross sequences, can be  
similarly generalized to our setting. 
We hope to explore possible roles played by Baxter algebras 
in other generalizations of the umbral calculus. 
We also plan to give a formulation of the umbral 
calculus in terms of coalgebras with operators 
by combining the approach 
in this paper and the coalgebra approach in \cite{RR,NS}. 
These projects will be carried out in subsequent papers. 

The layout of this paper is as follows. 
In section~\ref{sec:umb}, we review the umbral calculus 
and give a characterization of umbral calculus in terms of 
Baxter algebra. 
In section~\ref{sec:lambda}, we define sequences of 
$\lambda$-binomial type and formulate the basic theory of 
the $\lambda$-umbral calculus that generalizes the classic 
theory of Rota. 
In section~\ref{sec:rel}, under the assumption that 
$C$ is a $\QQ$-algebra, we give an explicit construction 
of sequences of $\lambda$-binomial type. 
We also study the relation between the $\lambda$-binomial 
sequences and the classic binomial sequences. 

\section{The umbral calculus and umbral algebra}
\mlabel{sec:umb}

\subsection{Background on the umbral calculus}
We first recall some background on the umbral calculus. 
See \cite{RKO,Rom} for more details. 

Let $\NN$ be the set of non-negative integers. 
Let $C[x]$ be the $C$-algebra of polynomials with coefficients 
in $C$. 
The main objects to study in the classical umbral calculus 
are special sequences of polynomials called sequences 
of binomial type and Sheffer sequences. 
\begin{defn}
\begin{enumerate}
\item
A sequence $\{p_n (x)\mid n\in\NN \}$ of polynomials in 
$C[x]$ is called {\bf a sequence of binomial type} if 
\[ p_n(x+y)= \sum_{k=0}^n p_{n-k}(x) p_k(y), 
\forall y\in C, n\in \NN.
\]
\item
Given a sequence of binomial type $\{p_n(x)\}$, 
a sequence $\{s_n (x)\mid n\in\NN \}$ of polynomials in 
$C[x]$ is called {\bf a Sheffer sequence} relative to 
$\{p_n(x)\}$ if 
\[ s_n(x+y)= \sum_{k=0}^n p_{n-k}(x) s_k(y), 
\forall y\in C, n\in \NN.
\]
\end{enumerate}
\end{defn}

In order to describe these sequences, Rota and his 
collaborators studied the dual of the $C$-module $C[x]$ and 
endowed the dual with a $C$-algebra structure, called 
the umbral algebra. 
It can be defined as follows. 
Let $\{t_n \mid n\in \NN\}$ be a sequence of symbols. 
Let $\calf$ be the $C$-module 
$\prod_{n\in \NN} C t_n$, where the addition and 
scalar multiplication are defined componentwise. 
Define a multiplication on $\calf$ by assigning 
\begin{equation}
t_m t_n = \binc{m+n}{m} t_{m+n},\ m,\ n\in \NN. 
\mlabel{eq:div}
\end{equation}
This makes $\calf$ into a $C$-algebra, with $t_0$ 
being the identity. 
The $C$-algebra $\calf$, together with the basis
$\{t_n\}$ is called the {\bf umbral algebra}. 
When $C$ is a $\QQ$-algebra, it follows from 
Eq (\ref{eq:div}) that 
\[ t_n=\frac{t_1^n}{n!},\ n\in \NN\]
and so $\calf \cong C[[t_1]]$ as a $C$-algebra. 
But we want to emphasize the special basis $\{t_n\}$. 

One then identifies $\calf$ with the dual $C$-module of 
$C[x]$ by taking $\{t_n\}$ to be the dual basis of 
$\{x^n\}$. In other words, $t_k$ is defined by  
\[ t_k : C[x] \to C,\ x^n \mapsto \delta_{k,n},\ 
k,\ n\in \NN. \]
Rota and his collaborators removed the mystery of 
sequences of binomial type and Sheffer sequences by 
showing that they have a simple characterization in 
terms of the umbral algebra. 

Let $f_n,\ n\geq 0$ be a pseudo-basis of $\calf$. 
This means that $\{f_n,\ n\geq 0\}$ is linearly independent 
and generates $\calf$ as a topology $C$-module 
where the topology on $\calf$ is defined by the 
filtration 
\[ F^n=\left \{\sum_{k=1}^\infty c_n t_n 
    \mid c_n =0, n\leq k \right \}.\]
A pseudo-basis $f_n,\ n\geq 0$ is called a {\bf divided 
power pseudo-basis} if 
\[ f_m f_n = \binc{m+n}{m} f_{m+n},\ m,\ n\geq 0. \]
Much of the foundation for the umbral calculus can be 
summarized in the following theorem. 
\begin{theorem}
{\bf \cite{RR}} 
Let $C$ be a $\QQ$-algebra. 
\begin{enumerate}
\item
A polynomial sequence $\{p_n(x)\}$ is of binomial 
type if and only if it is the dual basis of a divided 
power pseudo-basis of $\calf$.
\item
Any divided power pseudo-basis of $C[[t]]$ is of 
the form 
$ f_n(x)=\frac{f^n(t)}{n!}$ for some $f\in C[[t]]$ with  
$\deg f=1$ (i.e., $\dis{f(t)=\sum_{k=1}^\infty c_k t^k,\ 
    c_1\neq 0}$). 
\end{enumerate}
\end{theorem}
Sheffer sequences can be similarly described. 

\subsection{Baxter algebras}
We will give a characterization of the umbral 
calculus in terms of Baxter algebras. 
For this purpose 
we recall some definitions and basic properties on 
Baxter algebras. 
For further details, see \cite{G-K1}. 

\begin{defn}
Let $A$ be a $C$-algebra and let $\lambda$ be in $C$. 
\begin{enumerate}
\item
A $C$-linear operator $P: A\to A$ is called a {\bf Baxter 
operator of weight $\lambda$} if 
\[ P(x)P(y) = P(xP(y))+P(yP(x))+\lambda P(xy),\ x,\ y\in A.\]
\item
The pair $(A,P)$, where $A$ is a $C$-algebra and $P$ is 
a Baxter operator on $A$ of weight $\lambda$, 
is called a {\bf Baxter $C$-algebra of weight $\lambda$}. 
\end{enumerate}
\end{defn}
We often suppress $\lambda$ from the notations when there  
is no danger of confusion.

\begin{defn}
Let $A$ be a $C$-algebra. 
A free Baxter algebra on $A$ of weight $\lambda$ is 
a weight $\lambda$ Baxter algebra $(F_A,P_A)$ together 
with a $C$-algebra morphism $j_A: A\to F_A$ such that, 
for any weight $\lambda$ Baxter algebra $(R,P)$ and 
any $C$-algebra morphism $f: A\to R$, there is a unique 
weight $\lambda$ Baxter algebra morphism 
$\tilde{f}: (F_A,P_A) \to (R,P)$ with 
$\tilde{f} \circ j_A = f$. 
\end{defn}

Free Baxter algebras were constructed in \cite{G-K1}, 
generalizing the work of Cartier ~\cite{Ca} and 
Rota~\cite{Ro2}. 
In the special case when $A=C$, we have 

\begin{theorem}
\mlabel{thm:free}
{\rm \cite{G-K1}}
Let $U_\lambda C$ be the direct product 
$\prod_{n\in \NN} C u_n$ of the rank one free $C$-modules 
$C u_n,\ n\in \NN$.  
\begin{enumerate}
\item
With the product defined by 
\[ u_m u_n = \sum_{k=0}^m \binc{m+n-k}{n} \binc{n}{k} 
    \lambda^k u_{m+n-k},\ m,\ n\in \NN, \]
$U_\lambda C$ becomes a $C$-algebra. 
\item
The operator 
\[ P_C: U_\lambda C\to U_\lambda C,\ u_n \mapsto u_{n+1},\ n\in \NN \]
is a Baxter operator or weight $\lambda$ on the 
$C$-algebra $U_\lambda C$. 
Further, the pair $(U_\lambda C,P_C)$ is the free complete Baxter 
algebra on $C$. 
\end{enumerate}
\end{theorem}

\begin{prop}
\mlabel{pp:basis}
Fix $\lambda\in C$. 
Let $\{v_n\}_n$ be a pseudo-basis of $U_\lambda C$. 
the following statements are equivalent. 
\begin{enumerate}
\item
\[ v_m v_n = \sum_{k=0}^m \binc{m+n-k}{n} \binc{n}{k} 
    \lambda^k v_{m+n-k},\ m,\ n\in \NN. \]
\item
The operator 
\[ P_v: U_\lambda C\to U_\lambda C,\ v_n \mapsto v_{n+1},\ n\in \NN \]
is a Baxter operator or weight $\lambda$ on the 
$C$-algebra $U_\lambda C$. 
Further, the pair $(U_\lambda C,P_v)$ is the free complete Baxter 
algebra on $C$. 
\end{enumerate}
\end{prop}
\proofbegin
The proposition follows immediately from 
the definition. 
\proofend

\begin{defn}
Let $\{v_n\}_n$ be a pseudo-basis of $U_\lambda C$. 
If any of the equivalent conditions in 
Proposition~\ref{pp:basis} is true, we call  
$\{v_n\}_n$ a {\bf $\lambda$-divided power pseudo-basis} 
of $U_\lambda C$. 
\end{defn}

We can now give a characterization of the umbral algebra 
and polynomial sequences of binomial type. 

\begin{theorem}
\mlabel{thm:uc0}
\begin{enumerate}
\item
The umbral algebra is the free Baxter algebra of weight 
zero on $C$. 
\item 
Let $\{p_n(x)\}_{n\in\NN}$ be a sequence of polynomials 
in $C[x]$. 
Then $\{p_n(x)\}_n$ is of binomial type if and only if 
it is the dual basis of a divided pseudo-basis of $U_\lambda C$. 
\end{enumerate}
\end{theorem}

\section{The $\lambda$-umbral calculus}
\mlabel{sec:lambda}

We now develop the theory of the $\lambda$-umbral calculus.  
We construct the $\lambda$-umbral algebra and 
establish the relation between the $\lambda$-umbral calculus
and $\lambda$-umbral algebra. 
In the special case when $\lambda=0$, we have the 
theory started by Rota on the umbral calculus. 

\subsection{Definitions}
In view of Theorem~\ref{thm:uc0}, 
we will study the following generalization of 
the umbral calculus. 
In order to get interesting examples, we will work with 
$C[[x]]$ instead of $C[x]$. 

\begin{defn}
A sequence $\{p_n (x)\mid n\in\NN \}$ of power series in 
$C[[x]]$ is a sequence of {\bf $\lambda$-binomial type}
if 
\[ p_n(x+y)= \sum_{k=0}^n \lambda^k 
    \sum_{i=0}^n \binc{n}{i}\binc{i}{k} 
    p_{n+k-i}(x) p_i(y), 
\forall y\in C, n\in \NN.
\]
\end{defn}
When $\lambda=0$, we recover the sequences of binomial 
type. 
We also note the following symmetric property. 
\begin{lemma}
\mlabel{lem:sym}
A sequence $\{p_n (x)\mid n\in\NN \}$ of power series in 
$C[[x]]$ is a sequence of {$\lambda$-binomial type}
if and only if 
\[ p_n(x+y)= \sum_{k=0}^n \lambda^k 
    \sum_{i=0}^n \binc{n}{i}\binc{i}{k} 
    p_{i}(x) p_{p+k-i}(y), 
\forall y\in C, n\in \NN.
\]
\end{lemma}
\proofbegin
We have 
\allowdisplaybreaks{
\begin{eqnarray*}
\lefteqn{ \sum_{k=0}^n\lambda^k\sum_{i=0}^n \binc{n}{i}
    \binc{i}{k} p_{n+k-i}(x) p_i(y) 
= \sum_{i=0}^n \sum_{k=0}^n \lambda^k \binc{n}{i}
    \binc{i}{k} p_{n+k-i}(x) p_i(y)} \\
&=& \sum_{i=0}^n \sum_{k=0}^i \lambda^k \binc{n}{i}
    \binc{i}{k} p_{n+k-i}(x) p_i(y) \\
&=& \sum_{i=0}^n \sum_{j=n-i}^n \lambda^{j-n+i} \binc{n}{i}
    \binc{i}{j-n+i} p_{j}(x) p_i(y) \\
&& \ \ \ ({\rm use\ substitution\ } j=n+k-i,\ k=j-n+i) \\
&=& \sum_{i=0}^n \sum_{j=0}^n \lambda^{j-n+i} \binc{n}{i}
    \binc{i}{j-n+i} p_{j}(x) p_{i}(y) \\
&=& \sum_{j=0}^n \sum_{i=0}^n \lambda^{j-n+i} \binc{n}{i}
    \binc{i}{j-n+i} p_{j}(x) p_{i}(y) \\
&=& \sum_{j=0}^n \sum_{\ell=j-n}^j \lambda^\ell \binc{n}{i}
    \binc{i}{\ell} p_{j}(x) p_{n+\ell-j}(y) \\
&& \ \ \ \ ({\rm using\ substitution\ } 
    \ell=j-n+i,\ i=n+\ell-j)\\
&=& \sum_{j=0}^n \sum_{\ell=0}^n \lambda^\ell \binc{n}{i}
    \binc{i}{\ell} p_{j}(x) p_{n+\ell-j}(y) \\
&=& \sum_{\ell=0}^n \sum_{j=0}^n \lambda^\ell \binc{n}{i}
    \binc{i}{\ell} p_{j}(x) p_{n+\ell-j}(y). 
\end{eqnarray*}}
This proves the lemma. 
\proofend

\begin{defn}
Fix a $\lambda\in C$. 
The algebra $U_\lambda C =\prod_n C u_n$, together with 
the $\lambda$-divided power pseudo-basis 
$\{u_n\}_n$, is called the {\bf $\lambda$-umbral algebra}. 
\end{defn}

Fix a pseudo-basis $\frakq=\{q_n(x)\}$ of $C[[x]]$ 
of $\lambda$-binomial type (see \S \ref{sec:rel} for 
the existence of of such sequences). 
Let $\Cq$ be the submodule of $C[[x]]$ generated 
by the pseudo-basis $\frakq$. 
As in the 
case of $\lambda=0$, we identify $U_\lambda C$ with the dual $C$-module of 
$\Cq$ by taking $\{u_n\}$ to be the dual basis of $\{q_n(x)\}$. 
Thus each element of $U_\lambda C$ can be regarded as 
a functional on $\Cq$. We denote 
$<\ \mid\ >_\lambda$ for the resulting pairing 
\[ U_\lambda C \otimes \Cq\ \to C. \]
It is characterized by 
\begin{equation}
\lpair{u_n}{q_k} = \delta_{k,n}\ \forall\  
    n,\ k\in \NN. 
\mlabel{eq:pair}
\end{equation}

\begin{lemma}
\mlabel{lm:pair}
The paring $\lpair{}{}$ is pseudo-perfect in the sense that  
\begin{enumerate}
\item
for a fixed $u\in U_\lambda C$, if  
$ \lpair{u}{f(x)}=0,\ \forall f(x)\in \Cq $, 
then $u=0$, 
and 
\item
for a fixed $f(x)\in \Cq$, 
if $\lpair{u}{f(x)}=0,\ \forall u \in U_\lambda C $, 
then $ f(x)=0$. 
\end{enumerate}
\end{lemma}
\proofbegin
It follows from Eq.(\ref{eq:pair}) and the fact 
that $\{u_n\}$ and $\{q_k\}$ are pseudo-bases. 
\proofend

\subsection{Basic properties}

The following elementary properties can be proved in 
the same way as in the classical case. 

\begin{prop} 
\mlabel{pp:exp}
Let $\{p_n(x)\}\subset \Cq$ be a sequence of 
$\lambda$-binomial type. Let $\{v_n\}\subset U_\lambda C$ 
be the dual basis of $\{p_n(x)\}$. 
\begin{enumerate}
\item
{\bf (The Expansion Theorem)} 
For any $u\in U_\lambda C$, 
\[ u=\sum_{n=0}^\infty <u\mid p_n(x)> v_n. \]
\item 
{\bf (The Polynomial Expansion Theorem)}
For any $p(x)\in \Cq$, we have 
\[ p(x) = \sum_{n=0}^\infty <v_n\mid p(x)> p_n(x). \]
\end{enumerate}
\end{prop}

Note that $U_\lambda C$ acts on itself on the right, making $U_\lambda C$ 
a right $U_\lambda C$-module. This $U_\lambda C$-module structure, through 
the pairing $<\ |\  >_\lambda$, makes $\Cq$ into a left 
$U_\lambda C$-module. More precisely, for $u\in U_\lambda C$ and $f\in 
\Cq$, the element $uf\in \Cq$ is characterized by 
\begin{equation}
 < v\mid uf>_\lambda =<vu \mid f>_\lambda,\ \forall\  
    v\in U_\lambda C.
\mlabel{eq:act}
\end{equation}

\begin{prop}
\mlabel{pp:act}
With the notation above, we have 
\[ u_k q_n(x) = \sum_{i=0}^k \lambda^i 
    \binc{n}{k} \binc{k}{i} q_{n-k+i}(x), 
\forall\ k,\ n\in \NN.\]
\end{prop}

\proofbegin
For each $m\in\NN$, we have 
\begin{eqnarray*}
\lefteqn{ \lpair{u_m}{u_ks_n}=\lpair{u_mu_k}{s_n} }\\
&=& \lpair{\sum_{i=0}^m\lambda^i \binc{m+k-i}{m} 
    \binc{m}{i}u_{m+k-i}}{s_n} \\
&=& \sum_{i=0}^m \lambda^i 
    \binc{m+k-i}{m}\binc{m}{i} \delta_{m+k-i,n}\\
&=& \sum_{i=0}^m \lambda^i 
    \binc{m+k-i}{m}\binc{m}{i} \delta_{m,n-k+i}\\
&=& \sum_{i=0}^\infty \lambda^i
    \binc{m+k-i}{m}\binc{m}{i} \delta_{m,n-k+i}\ 
\ \ \ (i>m\Rightarrow \binc{m}{i} =0)  \\
&=& \sum_{i=0}^\infty \lambda^i
\binc{n}{n-k+i}\binc{n-k+i}{i} \delta_{m,n-k+i}\\
&& \ \ \ \ 
{\rm (the\ definition\ of\ \delta_{m,n-k+i} \Rightarrow 
    m=n-k+i)}\\
&=& \sum_{i=0}^\infty \lambda^i
\binc{n}{n-k+i}\binc{n-k+i}{i} \lpair{u_m}{s_{n-k+i}} 
\ \ \ ({\rm Eq.} (\ref{eq:pair}))\\
\end{eqnarray*}
\begin{eqnarray*}
&=& \sum_{i=0}^\infty \lambda^i
\binc{n}{k}\binc{k}{i} \lpair{u_m}{s_{n-k+i}} 
\ \ \ \ (\binc{n}{n-k+i}\binc{n-k+i}{i}
        =\binc{n}{k}\binc{k}{i})\\
&=& \sum_{i=0}^k \lambda^i
\binc{n}{k}\binc{k}{i} \lpair{u_m}{s_{n-k+i}} 
\ \ \ \ \ (i>k \Rightarrow \binc{k}{i}=0) \\
&=& \lpair{u_m}{\sum_{i=0}^k \lambda^i
\binc{n}{k}\binc{k}{i} s_{n-k+i}}. 
\end{eqnarray*}
Now the proposition follows from Lemma~\ref{lm:pair}.
\proofend
\begin{prop}
\mlabel{pp:prod1}
Let $\frakq=\{q_n(x)\}$ be the fixed sequence of $\lambda$-binomial 
type. For any $u,\ v\in U_\lambda C$, we have 
\begin{equation}
\lpair{uv}{q_n(x)}=\sum_{k=0}^n\sum_{i=0}^n 
    \lambda^k \binc{n}{i} \binc{i}{k} 
    \lpair{u}{q_{n+k-i}(x)} \lpair{v}{q_i(x)}.
\end{equation}
\end{prop}
\proofbegin
By Proposition~\ref{pp:exp}, we have 
\[ u=\sum_{n=0}^\infty \lpair{u}{q_n(x)} u_n,\]
and 
\[ v=\sum_{n=0}^\infty \lpair{v}{q_n(x)} u_n.\]
Then we have 
\allowdisplaybreaks{
\begin{eqnarray*}
\lefteqn{uv=\sum_{m=0}^\infty \sum_{k=0}^\infty 
    \lpair{u}{q_k(x)} \lpair{v}{q_{m-k}(x)} u_k u_{m-k}}\\
&=& \sum_{m=0}^\infty \sum_{k=0}^\infty 
    \lpair{u}{q_k(x)} \lpair{v}{q_{m-k}(x)} 
    \sum_{i=0}^k \lambda^i \binc{m-i}{k} \binc{k}{i} u_{m-i}\\
&& \ \ \ \ ({\rm Proposition~\ref{pp:basis}})\\
&=& \sum_{m=0}^\infty \sum_{k=0}^\infty 
    \lpair{u}{q_k(x)} \lpair{v}{q_{m-k}(x)} 
    \sum_{j=m-k}^m \lambda^{m-j} \binc{j}{k} \binc{k}{m-j} 
    u_j \\
&&\ \ \ \ {\rm (replacing\ } i {\rm\ by\ } m-j)\\
&=& \sum_{m=0}^\infty \sum_{j=0}^\infty 
    \sum_{k=m-j}^m
    \lpair{u}{q_k(x)} \lpair{v}{q_{m-k}(x)} 
    \lambda^{m-j} \binc{j}{k} \binc{k}{m-j} u_j \\
&&\ \ \ \ ({\rm exchanging\ the\ second\ and\ the\ third\ 
    sum})\\
&=& \sum_{j=0}^\infty \left (\sum_{m=j}^\infty 
    \sum_{k=m-j}^m
    \lpair{u}{q_k(x)} \lpair{v}{q_{m-k}(x)} 
    \lambda^{m-j} \binc{j}{k} \binc{k}{m-j}\right ) u_j \\
&&\ \ \ \ ({\rm exchanging\ the\ first\ and\ the\ second\ 
    sum}).
\end{eqnarray*}}
Since $\binc{k}{m-j}=0$ for $m>k+j$ and 
$\binc{j}{k}=0$ for $k>j$, we have 
\allowdisplaybreaks{
\begin{eqnarray*}
\lefteqn{\sum_{m=j}^\infty 
    \sum_{k=m-j}^m
    \lpair{u}{q_k(x)} \lpair{v}{q_{m-k}(x)} 
    \lambda^{m-j} \binc{j}{k} \binc{k}{m-j} }\\
&=& \sum_{m=j}^{2j}
    \sum_{k=m-j}^m
    \lambda^{m-j} \binc{j}{k} \binc{k}{m-j} 
    \lpair{u}{q_k(x)} \lpair{v}{q_{m-k}(x)} \\
&=& \sum_{t=0}^j \sum_{k=t}^{t+j} 
    \lambda^{t} \binc{j}{k} \binc{k}{t} 
    \lpair{u}{q_k(x)} \lpair{v}{q_{t+j-k}(x)} \\
&& ({\rm replacing\ } m {\rm\ by\ } t+j)\\
&=& \sum_{t=0}^j \sum_{k=0}^{j} 
    \lambda^{t} \binc{j}{k} \binc{k}{t} 
    \lpair{u}{q_k(x)} \lpair{v}{q_{t+j-k}(x)} \\
&& (k<t \Rightarrow \binc{k}{t}=0 {\rm\ and\ } 
k>j \Rightarrow \binc{j}{k}=0).
\end{eqnarray*}}
Thus 
\[
uv = \sum_{j=0}^\infty 
\left ( \sum_{t=0}^j \sum_{k=0}^{j} 
    \lambda^{t} \binc{j}{k} \binc{k}{t} 
    \lpair{u}{q_k(x)} \lpair{v}{q_{t+j-k}(x)}\right ) u_j. 
\]
Applying this to $q_n(x)$ and using Eq. (\ref{eq:pair}), 
we have 
\[ \lpair{uv}{q_n(x)} = 
\sum_{t=0}^n \sum_{k=0}^{n} 
    \lambda^{t} \binc{n}{k} \binc{k}{t} 
    \lpair{u}{q_k(x)} \lpair{v}{q_{n+t-k}(x)}. \]
Exchanging $u$ and $v$ in the equation, we have 
\begin{eqnarray*}
\lefteqn{ \lpair{uv}{q_n(x)}= \lpair{vu}{q_n(x)}}\\
&=& \sum_{t=0}^n \sum_{k=0}^{n} 
    \lambda^{t} \binc{n}{k} \binc{k}{t} 
    \lpair{v}{q_k(x)} \lpair{u}{q_{n+t-k}(x)} \\
&=& \sum_{t=0}^n \sum_{k=0}^{n} 
    \lambda^{t} \binc{n}{k} \binc{k}{t} 
    \lpair{u}{q_{n+t-k}(x)} \lpair{v}{q_k(x)}.
\end{eqnarray*}
This proves the proposition. 
\proofend
\begin{prop}
\mlabel{pp:prod2}
Let $\{p_n(x)\}\subset \Cq$ be a sequence of 
$\lambda$-binomial type and 
let $u$ and $v$ be in $U_\lambda C$. Then 
\[ \lpair{uv}{p_n(x)}=\sum_{i=0}^n \sum_{j=0}^{n} 
    \lambda^{i} \binc{n}{j} \binc{j}{i} 
    \lpair{u}{p_{n+i-j}(x)} \lpair{v}{p_j(x)}.
\]
\end{prop}
\proofbegin
We follow the case when $\lambda=0$ \cite{RR}. 
Let $C[[x,y]]$ be the $C$-module of power series in the 
variables $x$ and $y$. Since 
$ C[[x,y]]\cong C[[x]] \otimes_C C[[y]] $
and $q_n(x)$ is a pseudo-basis of $C[[x]]$, 
elements of the form $q_i(x)q_j(y),\ i,\ j\in \NN$ form 
a pseudo-basis of $C[[x,y]]$. 
Thus any element $p(x,y)$ of $C[[x,y]]$ can be expressed 
uniquely in the form 
\[ p(x,y) = \sum_{i,j} c_{i,j} q_i(x) q_j(y),\ c_{i,j}\in C.\]
For $u\in U_\lambda C$, define 
\[ u_x p(x,y) = \sum_{i,j} c_{i,j} \lpair{u}{q_i(x)} q_j(y)\]
and 
\[ u_y p(x,y) = \sum_{i,j} c_{i,j} q_i(x) \lpair{u}{q_j(y)}.\]
Then by Proposition~\ref{pp:prod1} and the fact that 
$q_n(x)$ is a $\lambda$-binomial sequence, we have 
\begin{eqnarray*}
 \lpair{uv}{q_n(x)}&=&   
\sum_{i=0}^n \sum_{j=0}^{n} 
    \lambda^{i} \binc{n}{j} \binc{j}{i} 
    \lpair{u}{q_{n+i-j}(x)} \lpair{v}{q_j(x)} \\
&=& u_x v_y q_n(x+y).
\end{eqnarray*}
Since $q_n(x)$ is a basis of $\Cq$, by 
the $C$-linearity of the maps 
$\lpair{uv}{},\ u_x$ and $v_y$, 
we have 
\[ \lpair{uv}{p(x)} = u_x v_y p(x+y) \]
for any $p(x)\in \Cq$. 
Since $p_n(x)$ is of $\lambda$-binomial type, we obtain 
\begin{eqnarray*}
\lefteqn{ \lpair{uv}{p_n(x)} = u_x v_y p_n(x+y) }\\
&=& u_x v_y \left ( \sum_{i=0}^n\sum_{j=0}^n \lambda^i 
    \binc{n}{j}\binc{j}{i} p_{n+i-j}(x) p_j(y) \right )\\
&=& \sum_{i=0}^n\sum_{j=0}^n \lambda^i 
    \binc{n}{j}\binc{j}{i} \lpair{u}{p_{n+i-j}(x)} 
    \lpair{v}{p_j(x)} .
\end{eqnarray*}
\proofend 

Let $c\in C$. Define the {\bf shift operator} $E^c$ 
on $\Cq$ by 
\[ (E^c f)(x)= f(x+c). \] 
An operator $L$ on $\Cq$ is called {\bf shift invariant} 
if 
\[ E^c L = L E^c,\ \forall\ c\in C.\]
\begin{prop}
\mlabel{pp:inv}
The elements in $U_\lambda C$, regarded as operators on 
$\Cq$, 
are shift invariant. 
\end{prop}
\proofbegin
We only need to check that $u_k$ are shift invariant when applied 
to $q_n(x)$. 
Then since $q_n(x)$ is a basis of $\Cq$, 
$u_k$ is shift invariant on $\Cq$. 
Since $u_k$ is a pseudo-basis of $U_\lambda C$, 
we see further that 
every element in $U_\lambda C$ is shift invariant on 
$\Cq$, as is needed. 

Fix $k,n\geq 0$. We have
\begin{eqnarray*}
\lefteqn{E^c u_k q_n(x) = 
    \sum_{s=0}^k\lambda^s \binc{n}{k} \binc{k}{s} 
    q_{n-k+s}(n+c)}\\
&=& \sum_{s=0}^k\lambda^s \binc{n}{k} \binc{k}{s} 
    \sum_{\ell=0}^{n-k+s} \binc{n-k+s}{i}\binc{i}{\ell} 
    q_{n-k+s+\ell-i} q_i(c)\\
&=& \sum_{i=0}^n \sum_{s=0}^k \sum_{\ell=0}^n 
    \lambda^{s+\ell} \binc{n}{k}\binc{k}{s} 
    \binc{n-k+s}{i}\binc{i}{\ell} q_{n+\ell+s-k-i}q_i(c)\\
&=& \sum_{i=0}^n \sum_{w=0}^{n+k} \lambda^w\binc{n}{k}
    \sum_{s=0}^w \binc{k}{s} \binc{n-k+s}{i} \binc{i}{w-s} 
    q_{n+w-k-i}(x) q_i(c)\\
&& ({\rm letting\ } w=\ell+s).
\end{eqnarray*}
Also
\allowdisplaybreaks{
\begin{eqnarray*}
\lefteqn{u_k E^c q_n(x)=u_k q_n(x+c)}\\
&=& u_k \sum_{\ell=0}^n \lambda^\ell \sum_{i=0}^n 
 \binc{n}{i} \binc{i}{\ell} q_{n+\ell-i}(x) q_i(c)\\
&=& \sum_{\ell=0}^n\lambda^\ell\sum_{i=0}^n
\binc{n}{i}\binc{i}{\ell} \sum_{s=0}^k\lambda^s 
    \binc{n+\ell-i}{k}\binc{k}{s} q_{n+\ell-i-k+s} q_i(c)\\
&=& \sum_{i=0}^n\sum_{s=0}^k\sum_{\ell=0}^n 
    \lambda^{\ell+s}\binc{n}{i}\binc{i}{\ell} 
    \binc{n+\ell-i}{k}\binc{k}{s} q_{n+\ell-i-k+s}q_i(c)\\
&=& \sum_{i=0}^n\sum_{w=0}^{n+k}\lambda^w\binc{n}{k} 
    \sum_{s=0}^w\binc{n+w-s-i}{w-s}\binc{n-k}{i-w+s} 
    \binc{k}{s} q_{n+w-i-k}q_i(c)\\
&& ({\rm letting\ } w=\ell+s)\\
&=& \sum_{i=0}^n\sum_{w=0}^{n+k}\lambda^w\binc{n}{k} 
    \sum_{s=0}^w\binc{k}{w-s}\binc{n+s-i}{s} 
    \binc{n-k}{i-s} q_{n+w-i-k}q_i(c)\\
&& \ \ \ \ ({\rm replacing\ } s {\rm\ by\ } w-s).
\end{eqnarray*}}
Thus we only need to prove 
\begin{equation}
 \sum_{s=0}^w \binc{k}{s} \binc{n-k+s}{i} \binc{i}{w-s}
=\sum_{s=0}^w \binc{k}{w-s}\binc{n+s-i}{s} \binc{n-k}{i-s} 
\mlabel{eq:inv}
\end{equation}
for $n,k,i,w\geq 0$. 
With the help of the Zeilberger algorithm\cite{PWZ}, 
we can verify that both sides of the equation satisfy 
the same recursive relation 
\begin{eqnarray*}
&&((k+i-w)*(k-n+i-w-1)\\
&&+(k^2-k*n+k*i-3*k*w
-n*i+2*n*w+i^2-3*i*w\\
&&+2*w^2-4*k+2*n-4*i+5*w+3)*W\\
&&-(w+2)*(k-n+i-w-2)*W^2)F(n,k,i,w)=0
\end{eqnarray*}
for $n,k,i\geq 0$. 
Here $W$ is the shift operator $WF(n,k,i,w)=F(n,k,i,w+1)$. 
we can also easily verify Eq.(\ref{eq:inv}) directly 
for $w=0,1$. Therefore 
Eq. (\ref{eq:inv}) is verified. This completes the 
proof of Proposition~\ref{pp:inv}.
\proofend

\subsection{Sequences of $\lambda$-binomial type and 
    Baxter bases}
Now we are ready to state our main theorem in the theory 
of the $\lambda$-umbral calculus. 

\begin{theorem}
\mlabel{thm:ucl}
Let $\{p_n(x)\}_{n\in \NN}$ be a basis of 
$\Cq$. The following statements are equivalent. 
\begin{enumerate}
\item
The sequence $\{p_n (x)\}$ is a sequence of $\lambda$-binomial type. 
\item
The sequence $\{p_n(x)\}$  is the dual basis of 
a $\lambda$-divided power pseudo-basis of $U_\lambda C$. 
\item
For all $u$ and $v$ in $U_\lambda C$, 
\[ \lpair{uv}{p_n(x)}=\sum_{i=0}^n\sum_{j=0}^n 
    \lambda^i \binc{n}{j} \binc{j}{i} 
    \lpair{u}{p_{n+i-j}(x)} \lpair{v}{p_j(x)}.
\]
\end{enumerate}
\end{theorem}
\begin{remark}
As in the case when $\lambda=0$, the third statement 
in the theorem 
has the following interpretation in terms of coalgebra. 

\begin{enumerate}
\item[3']
The $C$-linear map $q_n(x)\mapsto p_n(x),\ n\in \NN,$ defines 
an automorphism of the $C$-coalgebra $\Cq$. 
Here the coproduct 
\[ \Delta: \Cq \to \Cq \otimes \Cq \]
is defined by first assigning 
\[ \Delta(q_n(x))= \sum_{i=0}^n\sum_{j=0}^n 
    \lambda^i \binc{n}{j} \binc{j}{i} 
    q_{n+i-j}(x)\otimes q_j(x),\ \ n\in \NN\]
and then extend $\Delta$ to $\Cq$ by $C$-linearity. 
\end{enumerate}
\end{remark}

\proofbegin
($1\Rightarrow 3$) has been proved in 
Proposition~\ref{pp:prod2}. 

($3\Rightarrow 2$): 
Let $v_n$ be the dual basis of $p_n(x)$ in $U_\lambda C$. 
Then we have 

\begin{eqnarray*}
\lpair{v_mv_n}{p_k(x)}&=&\sum_{i=0}^k\sum_{j=0}^k 
    \lambda^i \binc{k}{j} \binc{j}{i} 
    \lpair{v_m}{p_{k+i-j}(x)} \lpair{v_n}{p_j(x)}\\
&=& \sum_{i=0}^k\sum_{j=0}^k 
    \lambda^i \binc{k}{j} \binc{j}{i} 
    \delta_{m,k+i-j} \delta_{n,j}\\
&=& \sum_{i=0}^k 
    \lambda^i \binc{k}{m} \binc{m}{i} 
    \delta_{m,k+i-n}\ \ \ 
(\delta_{n,j}=0 {\rm\ for\ } j\neq m)\\
&=& \sum_{i=0}^k 
    \lambda^i \binc{k}{m} \binc{m}{i} 
    \delta_{m+n-k,i}\\
&=& \lambda^{m+n-k}\binc{k}{m} \binc{m}{m+n-k}.
\end{eqnarray*}
On the other hand,  
\begin{eqnarray*}
\lefteqn{ \lpair{\sum_{i=0}^m \lambda^i \binc{m+n-i}{m}
    \binc{m}{i} v_{m+n-i}}{p_k(x)} }\\
&=& \sum_{i=0}^m \lambda^i \binc{m+n-i}{m}
    \binc{m}{i} \delta_{m+n-i,k}\\
&=& \sum_{i=0}^m \lambda^i \binc{m+n-i}{m}
    \binc{m}{i} \delta_{m+n-k,i}\\
&=& \lambda^{m+n-k} \binc{k}{m}
    \binc{m}{m+n-k}. 
\end{eqnarray*}
So by Lemma~\ref{lm:pair}, 
\[ v_m v_n =\sum_{i=0}^m \lambda^i \binc{m+n-i}{m}
    \binc{m}{i} v_{m+n-i} \]
and $v_n$ is a $\lambda$-divided power pseudo-basis of $U_\lambda C$. 

($2\Rightarrow 1$) 
Let $p_n(x)$ be the dual basis of a $\lambda$-divided 
pseudo-basis $v_n$ 
of $U_\lambda C$. Let $y$ be a variable and consider the 
$C$-algebra $C_1\defeq C[[y]]$. Regarding $p_n(x)$ as 
elements in $C_1[[x]]$, the sequence $p_n(x)$ is still of 
$\lambda$-binomial type. 
Also, let 
\[ U_\lambda C_1= \prod_{n\in\NN} C_1 u_n\]
be the $C_1$-algebra obtained from $U_\lambda C$ by scalar extension. Then 
$v_n$ is also a $\lambda$-divided power pseudo-basis of $U_\lambda C_1$. Further, the 
$C$-linear perfect pairing (\ref{eq:pair}) extends to a 
$C_1$-linear perfect pairing between $C_1\!<\!\frakq\!>$ 
and $U_\lambda C_1$. Under 
this pairing, $\{p_n(x)\}$ is still 
the dual basis of $\{v_n\}$. Note 
that all previous results applies when $C$ is replaced by $C_1$ 
since we did not put any restrictions on $C$. With this in mind, 
we have 
\begin{eqnarray*}
\lefteqn{p_n(x+y)=\sum_{k=0}^\infty \lpair{u_k}{p_n(x+y)} 
    p_k(x)\ \ \ ({\rm Lemma~\ref{pp:exp}})} \\
&=& \sum_{k=0}^\infty \lpair{u_k}{E^y p_n(x)} p_k(x)
    \ \ \ \ ({\rm definition\ of\ } E^y)\\
&=& \sum_{k=0}^\infty \lpair{u_kE^y }{p_n(x)} p_k(x)
    \ \ \ \ ({\rm Eq. (\ref{eq:act})})\\
&=& \sum_{k=0}^\infty \lpair{E^y u_k }{p_n(x)} p_k(x)
    \ \ \ \ ({\rm Proposition~\ref{pp:inv}})\\
&=& \sum_{k=0}^\infty \lpair{E^y}{u_k p_n(x)} p_k(x)
    \ \ \ \ ({\rm Eq.~\ref{eq:act}})\\
&=& \sum_{k=0}^\infty \lpair{E^y}
    {\sum_{i=0}^k \lambda^i 
    \binc{n}{k} \binc{k}{i} p_{n-k+i}(x)} p_k(x)
    \ \ \ \ ({\rm Proposition~\ref{pp:act}})\\
&=& \sum_{k=0}^\infty 
 \sum_{i=0}^k \lambda^i 
    \binc{n}{k} \binc{k}{i} p_{n-k+i}(y) p_k(x)
    \ \ \ \ ({\rm definition\ of\ }E^y)\\
&=& \sum_{k=0}^n 
 \sum_{i=0}^n \lambda^i 
    \binc{n}{k} \binc{k}{i} p_{n-k+i}(y) p_k(x)
    \ \ \ \ (k>n\Rightarrow \binc{n}{k}=0,\ 
    i>k\Rightarrow \binc{k}{i}=0)\\
&=& \sum_{i=0}^n 
 \sum_{k=0}^n \lambda^i 
    \binc{n}{k} \binc{k}{i} p_k(x)p_{n+i-k}(y) 
\ \ \ \ \ ({\rm exchanging\ the\ sums}).
\end{eqnarray*}
\proofend

\section{$\lambda$-binomial sequences in 
$C[[x]]$}
\mlabel{sec:rel}

In this section we study $\lambda$-binomial sequences 
in $\Cq$ when $C$ is a $\QQ$-algebra $C$. 
We show that $U_\lambda C$ is isomorphic to $C[[t]]$ 
regardless of $\lambda$. 
We further show that the pairing in 
the $\lambda$-umbral calculus and the pairing in the 
classical umbral calculus (when $\lambda=0$) are 
compatible. The situation is quite 
different when $C$ is not a $\QQ$-algebra and will 
be studied in a subsequent paper.

\subsection{The $\lambda$-umbral algebra and 
    $\lambda$-divided powers}
We first describe the $\lambda$-umbral algebra 
$U_\lambda C$. 

\begin{defn}
An power series $f(t)\in C[[t]]$ is call {\bf delta} 
if $\deg f=1$, that is, $f(t)=\sum_{k=1}^\infty a_k t^k$ 
with $a_1\neq 0$. 
\end{defn}

\begin{prop}
\mlabel{pp:power}
\begin{enumerate}
\item
Let $f(t)$ be a delta series. 
Then 
$$ \tau_n(f(t))\defeq \frac{f(t)(f(t)-\lambda)\cdots 
    (f(t)-\lambda(n-1))}{n!}, 
    n\geq 0$$
form a $\lambda$-divided power pseudo-basis for 
$C[[t]]$. 
\item
The map 
$$ C[[t]] \to U_\lambda C,\ 
\tau_n(t)=\frac{t(t-\lambda)\cdots(t-\lambda(n-1))}{n!} 
    \mapsto t_n,\ n\geq 0$$
identifies $C[[t]]$ with the $\lambda$-umbral algebra 
$U_\lambda C$. 
\end{enumerate}
\end{prop}
\proofbegin
1. Since $f(t)$ is delta, we have $\deg \tau_n(f)=n$. 
So $\tau_n(f(t))$ form a pseudo-basis for $C[[t]]$. 
Thus we only need to show that the pseudo-basis is 
a $\lambda$-divided power series. 
For this 
we prove by induction on $m\in\NN$ that, for any $n\in \NN$,
\begin{equation}
 \tau_m(f(t))\tau_n(f(t))= \sum_{k=0}^m \bincc{m+n-k}{m}\bincc{m}{k}
    \lambda^k \tau_{m+n-k}(f(t)).
    \mlabel{eq:div1}
\end{equation}
The equation clearly holds for $m=0$. 
Assume that equation~(\ref{eq:div1}) holds for an $m\in \NN$ 
and all $n\in \NN$. 
Then for any $n\in\NN$, we have
\allowdisplaybreaks{
\begin{eqnarray*}
\lefteqn{\sum_{k=0}^{m+1} \bincc{m+n+1-k}{m+1}\bincc{m+1}{k}
    \lambda^k \tau_{m+n+1-k}(f(t))}\\
&=& \sum_{k=0}^{m+1}
    \bincc{m+n+1-k}{m+1}
    \bincc{m}{k} \lambda^k \tau_{m+n+1-k}(f(t)) \\
&+& \sum_{k=0}^{m+1}
    \bincc{m+n+1-k}{m+1}
    \bincc{m}{k-1} \lambda^k \tau_{m+n+1-k}(f(t)) \\
&=& \sum_{k=0}^{m}
    \bincc{m+n+1-k}{m+1}
    \bincc{m+1}{k} \lambda^k \tau_{m+n+1-k}(f(t)) \\
&&+ \sum_{k=0}^{m}
    \bincc{m+n-k}{m+1}
    \bincc{m+1}{k} \lambda^{k+1} \tau_{m+n-k}(f(t)) \\
&=& \sum_{k=0}^{m}
    \bincc{m+n-k}{m} \bincc{m}{k}
    \lambda^k \tau_{m+n-k}(f(t)) (\frac{f(t)-(m+n-k)\lambda}{m+1})\\
&&+ \sum_{k=0}^{m}
    \bincc{m+n-k}{m} \bincc{m}{k} \lambda^k
    \tau_{m+n-k}(f(t)) (\frac{(n-k)\lambda}{m+1})\\
&=& \sum_{k=0}^{m}
    \bincc{m+n-k}{m} \bincc{m}{k}
    \lambda^k \tau_{m+n-k}(f(t))
(\frac{f(t)-(m+n-k)\lambda}{m+1}-\frac{(n-k)\lambda}{m+1})\\
&=& \tau_m(f(t))\tau_n(f(t))(\frac{f(t)-m\lambda}{m+1})\\
&=& \tau_{m+1}(f(t))\tau_n(f(t)). 
\end{eqnarray*}}
This completes the induction.

\noindent
2. The statement is clear since 
$\{\tau_{\lambda,n}(t)\}_n$ is a weight $\lambda$ 
$\lambda$-divided power pseudo-basis. 
\proofend

We next construct a sequence of $\lambda$-binomial type. 
For a given $\lambda\in C$, 
we use 
$\displaystyle{e_\lambda(x)=
    \frac{e^{\lambda x}-1}{\lambda}}$
to denote the series  
$\displaystyle{\sum_{k=1}^\infty \frac{\lambda^{k-1} x^k}{k!}.}$
When $\lambda=0$, we have 
$e_\lambda(x)=x$.

\begin{prop}
\mlabel{pp:ex}
Let $C$ be a $\QQ$-algebra. $\{e_\lambda^n(x)\}_n$
is a $\lambda$-binomial pseudo-basis for $C[[x]]$. 
\end{prop}

\proofbegin
We prove by induction the equation 
\begin{equation}
e_\lambda^n(x+c)= \sum_{k=0}^n \lambda^k 
    \sum_{i=0}^n \binc{n}{i}\binc{i}{k} 
    e_\lambda^{n+k-i}(x) e_\lambda^i(c), 
\forall c\in C, n\in \NN.
\mlabel{eq:bas}
\end{equation}
Clearly Eq. (\ref{eq:bas}) holds for $n=0$. 
It is also easy to verify the equation
$$e_\lambda(x+c)=e_\lambda(x)+e_\lambda(c)
    +\lambda e_\lambda(x)e_\lambda(c)$$
which is Eq. (\ref{eq:bas}) when $n=1$.
Using this and the induction hypothesis, we get 

\begin{eqnarray*}
\lefteqn{e_\lambda^{n+1}(x+c)=e(x+c)\sum_{k=0}^n\lambda^k 
    \sum_{i=0}^n \binc{n}{i}\binc{i}{k}e_\lambda^{n+k-i}(x)e_\lambda^i(c)}\\
&=& \sum_{k=0}^n\lambda^k \sum_{i=0}^n 
    \binc{n}{i}\!\!\binc{i}{k} 
    \lleft e_\lambda^{n+1+k-i}(x)e_\lambda^i(c) +e_\lambda^{n+k-i}(x)e_\lambda^{i+1}(c) 
    +\lambda e_\lambda^{n+1+k-i}e_\lambda^{i+1}(c)\lright  \\
&=& \sum_{k=0}^n\lambda^k \sum_{i=0}^n 
    \binc{n}{i}\!\!\binc{i}{k} 
    e_\lambda^{n+1+k-i}(x)e_\lambda^i(c) 
+ \sum_{k=0}^n \lambda^k \sum_{i=1}^{n+1}  
    \binc{n}{i-1}\!\!\binc{i-1}{k} e_\lambda^{n+1+k-i}(x)e_\lambda^i(c) \\
&& +\sum_{k=1}^{n+1} \lambda^k  \sum_{i=1}^{k+1}
    \binc{n}{i-1}\binc{i-1}{k} e_\lambda^{n+2+k-i}e_\lambda^i(c)\\ 
&=& \sum_{i=0}^{n+1} \left ( \binc{n}{i}+\binc{n}{i-1} \right ) 
    e_\lambda^{n+1-i}(x) e_\lambda^i(c) \\
&&+\sum_{k=1}^n \lambda^k  \lleft e_\lambda^{n+1}(x) 
    +\sum_{i=1}^n  \!\!
    \left (\!\! \binc{n}{i}\!\! \binc{i}{k}\!\! + \!\!
    \binc{n}{i-1}\!\!\binc{i-1}{k}\!\!
    +\!\!\binc{n}{i-1}\!\!\binc{i-1}{k-1} \!\!
    \right ) e_\lambda^{n+1+k-i}(x) e_\lambda^i (c) \\
&&+\left ( \binc{n}{k} +\binc{n}{k-1} \right ) 
     e_\lambda^k(x) e_\lambda^{n+1}(c) \lright  
+\lambda^{n+1} \sum_{i=0}^{n+1} \binc{n}{i-1}\binc{i-1}{n}
    e_\lambda^{2n+2-i}(x)e_\lambda^i(c)\\
&=& \sum_{i=0}^{n+1} \binc{n+1}{i} e_\lambda^{n+1-i}(x) e_\lambda^i(c) 
+\sum_{k=1}^n \lambda^k \sum_{i=0}^{n+1}\binc{n+1}{i} 
    \binc{i}{k} e_\lambda^{n+1+k-i}(x) e_\lambda^i(c) \\
&&+\lambda^{n+1}\sum_{i=0}^{n+1} \binc{n+1}{i}\binc{i}{n+1} 
    e_\lambda^{2n+2-i}(x) e_\lambda^i(c) \\
&=& \sum_{k=0}^{n+1} \lambda^{k} \sum_{i=0}^{n+1} 
    \binc{n+1}{i}\binc{i}{k} e_\lambda^{n+1+k-i}(x) e_\lambda^i(c).
\end{eqnarray*}
Thus Eq. (\ref{eq:bas}) is proved.
\proofend

Thus we can use the $\lambda$-binomial pseudo-basis 
$\frakq=\{e_\lambda^n(x)\}$
of $C[[x]]$ and the $\lambda$-divided power pseudo-basis 
$\displaystyle{\tau_n(t)=\frac{t(t-a)\cdots (t-a(n-1))}{n!},} 
\ n\geq 0,$
of $C[[t]]$ to define the pairing 
$$<\ |\ >_\lambda: U_\lambda C \times \Cq 
    \to C$$ 
as is described in Eq.(\ref{eq:pair}). 

\begin{defn}
Let $f(t)\in C[[t]]$ be a delta series. 
The dual basis of $\{\tau_n(f(t))\}_n$ in $\Cq$ is 
called the {\bf associated sequence for $f(t)$}. 
\end{defn}

We show that every $\lambda$-binomial sequence is 
associated to a delta power series, as in the case 
when $\lambda=0$. 

\begin{theorem}
\mlabel{thm:char0}
Let $C$ be a $\QQ$-algebra. 
Let $\{s_n(x)\}_n$ be a basis of $\Cq$. 
Then 
$\{s_n(x)\}$ is a sequence of $\lambda$-binomial type if 
and only if 
$\{s_n(x)\}$ is the associated sequence of a delta 
series $f(t)$ in $C[[t]]$. 
\end{theorem}

\proofbegin
If $\{s_n(x)\}$ is the associated sequence of a 
delta series $f(t)$, then by Theorem~\ref{thm:ucl}  
and Proposition~\ref{pp:power}, 
$\{s_n(x)\}$ is a $\lambda$-binomial sequence. 

Conversely, let $\{s_n(x)\}$ be a $\lambda$-binomial 
sequence and let $\{f_n(t)\}$ be the dual basis 
of $\{s_n(x)\}$ in $C[[t]]$. It follows from 
Theorem~\ref{thm:ucl} that $\{f_n(t)\}$ is a 
$\lambda$-divided power series. We only need to show 
\begin{enumerate}
\item
$f_1(t)$ is a delta series, and 
\item
$f_n(t)=\tau_n(f_1)(t), \forall n\geq 0.$
\end{enumerate}
We first prove (2) by induction. 
Since $\{f_n(t)\}$ is a $\lambda$-divided power basis, 
we have 
$$ f_m(t)f_n(t)= \sum_{k=0}^m \bincc{m+n-k}{m}\bincc{m}{k}
    \lambda^k f_{m+n-k}(t).
$$
Taking $m=n=0$, we have 
$f_0(t)f_0(t)=f_0(t)$. Thus $f_0(t)=1$. 
Assuming $f_n(t)=\tau_n(f_1)(t)$, from
$f_n(t)f_1(t)=f_{n+1}(t)+ n\lambda f_n(t)$
we have 
$$f_{n+1}(t)=f_n(t)(f_1(t)-n\lambda)
    =\tau_n(f_1)(t) (f_1(t)-n\lambda) 
    =\tau_{n+1}(f_1)(t).$$
This proves (2). 
Then (1) follows since if $f_1(t)$ is not a delta series, 
then $\{\tau_n(f)(t)\}$ cannot be a pseudo-basis 
of $C[[t]]$. 
\proofend

\subsection{Compatibility of $\lambda$-calculus}

We now show that, for any given $\lambda\in C$, 
the $\lambda$-umbral calculus is compatible with the 
classical umbral calculus (when $\lambda=0$). 
The compatibility is in the following sense. 
Since $\{x^n\}$ and $\frakq=\{q_n\}$ are 
pseudo-basis of $C[[x]]$, the pairings 
$$<\ |\ >_0: C[t]\otimes C[x] \to C$$ 
and 
$$<\ | \ >_\lambda: C[t]\otimes \Cq \to C$$
extend to pairings  
$$[\ \ |\ \ ]_0: C[t] \otimes C[[x]] \to C$$
and 
$$[\ \ |\ \ ]_\lambda: C[t] \otimes C[[x]] \to C.$$ 
\begin{defn}
We say that $<\ |\ >_0$ and $<\ |\ >_\lambda$ are 
compatible if 
$$ [\ \ |\ \ ]_0=[\ \ |\ \ ]_\lambda.$$
\end{defn}

\begin{theorem}
\mlabel{thm:comp}
The pairings 
$$<\ | \ >_0: C[t] \otimes C[[x]] \to C$$ 
and  
$$<\ | \ >_\lambda: C[t] \otimes C[[x]] \to C$$
are compatible. 
\end{theorem}
This theorem has the following immediate consequence. 

\begin{coro}
\mlabel{co:comp}
Let $C$ be a $\QQ$-algebra. A pseudo-basis $\{p_n(x)\}$ 
of $C[[x]]$ is of $\lambda$-binomial type if and only 
if it is the dual basis of a $\lambda$-divided power 
pseudo-basis of $C[[t]]$ under the pairing 
$<\ |\ >_0$.
\end{coro}

\noindent
{\bf Proof of Theorem~\ref{thm:comp}:}
Since $\{\tau_{\lambda,n}(t)\}_n$ is a basis of 
$C[t]$ and $\{e_\lambda^k(x)\}_k$ is a pseudo-basis of 
$C[[x]]$, we only need to prove 

\begin{equation}
[ \tau_{\lambda,n}(t)\ |\ e_\lambda^k(x)]_0 =\delta_{n,k}, \
n,k\geq 0. 
\mlabel{eq:comp}
\end{equation}
We will prove this equation by induction on $n$. 
When $n=0$, $\tau_{\lambda,n}(t)=1$. 
Since $\deg e_\lambda(x)=1$, we have 
$[\tau_{\lambda,0}(t)\ |\ e_\lambda^k(x)]_0=\delta_{0,k}$. 
Assume that Eq. (\ref{eq:comp}) holds for an $n\geq 0$ and 
all $k\geq 0$. 
Then we have 
\begin{eqnarray*}
\lefteqn{ [\tau_{\lambda,n+1}(t)\ |\ e_\lambda^k(x)]_0 }\\
&=& 
[\frac{t}{n+1}\tau_{\lambda,n}(t)-\frac{n\lambda}{n+1}
\tau_{\lambda,n}(t)\ |\ e_\lambda^k(x)]_0\\
&=&
\frac{1}{n+1}[t\tau_{\lambda,n}(t)\ |\ e_\lambda^k(x)]_0 
    -\frac{n\lambda}{n+1}[\tau_{\lambda,n}(t)\ |\ e_\lambda^k(x)]_0\\
&=& \frac{1}{n+1}[t\tau_{\lambda,n}(t)\ |\ e_\lambda^k(x)]_0 
    - \frac{n\lambda}{n+1} \delta_{n,k}.
\end{eqnarray*}
By \cite[Theorem 2.2.5]{Rom}, we have 
$$[t\tau_{\lambda,n}(t)\ |\ e_\lambda^k(x)]_0 
=[\tau_{\lambda,n}(t)\ |\ \frac{d}{dx} e_\lambda^k(x)]_0.$$ 
Since 
$$\frac{d}{dx} e_\lambda^k(x)=ke_\lambda^{k-1}(x)e^{\lambda x} 
=\lambda k e_\lambda^k(x)+ke_\lambda^{k-1}(x), $$
we have 
\begin{eqnarray*}
\lefteqn{[\tau_{\lambda,n+1}(t)\ |\ e_\lambda^k(x)]_0}\\
&=& \frac{1}{n+1}[\tau_{\lambda,n}(t)\ |\ \lambda k e_\lambda^k(x) 
    +k e_\lambda^{k-1}(x)]_0-\frac{n\lambda}{n+1}\delta_{n,k}\\
&=& \frac{\lambda k}{n+1} \delta_{n,k} 
+\frac{k}{n+1}\delta_{n,k-1} -\frac{n\lambda}{n+1} \delta_{n,k}\\
&=& 
\frac{\lambda(k-n)}{n+1}\delta_{n,k}+\frac{k}{n+1}\delta_{n,k-1}\\
&=& \delta_{n+1,k}.
\end{eqnarray*}
This completes the induction. 
\proofend

\noindent
{\bf Acknowledgement: }
The author thanks William Keigher for helpful discussions.

\addcontentsline{toc}{section}{\numberline {}References}

\end{document}